\documentclass[12pt]{amsart}
\usepackage{a4}
\usepackage{amssymb}
\usepackage{amscd}
\newtheorem{thm}{Theorem}
\newtheorem{fact}{Fact}

\newtheorem{lem}[thm]{Lemma}

\numberwithin{equation}{section}

\newcommand{\norm}[1]{\left\Vert#1\right\Vert}
\newcommand{\abs}[1]{\left\vert#1\right\vert}

\newcommand{\N}{\mathcal{N}}

\newcommand{\U}{\mathcal{U}}
\newcommand{\M}{\mathcal{M}}

\newcommand{\K}{\mathcal{K}}
\begin{document}
\title{Extremal cases of exactness constant and completely bounded projection constant}
\author{Hun Hee Lee}
\address{Department of Mathematical Sciences, Seoul National University
          San56-1 Shinrim-dong Kwanak-gu Seoul 151-747, Korea}
\email{bbking@amath.kaist.ac.kr; lee.hunhee@gmail.com}
\thanks{This research is supported by BK21 program}

\begin{abstract}
We investigate some extremal cases of exactness constant and completely bounded projection constant.
More precisely, for an $n$-dimensional operator space $E$ we prove that $\lambda_{cb}(E) = \sqrt{n}$
if and only if $ex(E) = \sqrt{n}$, which is equivalent to $\lambda_{cb}(E) < \sqrt{n}$ if and only if $ex(E) < \sqrt{n}$.

\end{abstract}
\maketitle

\section{Introduction}
Exactness constant and completely bounded (shortly c.b.) projection constant are fundamental
quantities in operator space theory. 

For an operator space $E \subseteq B(H)$, the {\it c.b. projection constant} of $E$,
$\lambda_{cb}(E)$ is defined by
$$\lambda_{cb}(E) = \inf\{ \norm{P}_{cb} \,|\, P :B(H) \rightarrow E, \,\,\text{projection onto}\,\, E \}.$$

Let $B = B(\ell_2)$ and $\K$ be the ideal of all compact operators on $\ell_2$, and let
$$T_E : (B\otimes_{\min}E) / (\K\otimes_{\min}E) \rightarrow (B / \K)\otimes_{\min}E$$
be the map obtained from $$q \otimes I_E : B\otimes_{\min}E \rightarrow (B / \K)\otimes_{\min}E$$
by the taking quotient with respect to $\K\otimes_{\min}E$, where $q : B \rightarrow \K$ is the canonical quotient map.
Then the {\it exactness constant} of $E$, $ex(E)$ is defined by $$ex(E) = \norm{T^{-1}_E}.$$
It is well known that the exactness constant is the same with $d_{\mathcal{S}\K}(E)$, where
$$d_{\mathcal{S}\K}(E) = \inf\{ d_{cb}(E,F) : F \subseteq \K \}$$ when $E$ is finite dimensional. (\cite{P0})

The followings are well known facts about these quantities
(Chapter 7 and 17 of \cite{P3} and section 9 of \cite{P1}):

\begin{fact}\label{fact1} For a finite dimensional operator space $E$ we have $$ex(E) = d_{\mathcal{S}\K}(E) \leq \lambda_{cb}(E).$$
\end{fact}

\begin{fact}\label{fact2} When $\text{dim}(E) = n \in \mathbb{N}$, we have $$\lambda_{cb}(E) \leq \sqrt{n}.$$
\end{fact}

Thus we can say that, for an $n$-dimensional operator space $E$, $\lambda_{cb}(E)$ and $ex(E)$ are both bounded by $\sqrt{n}$,
and this upper bound is known to be asymptotically sharp. Indeed, we have
$ex(\max \ell^n_1) \geq \frac{n}{2\sqrt{n-1}} \,\, \text{for} \,\, n \geq 2$. (\cite{P0})
However, we do not know whether there is an
$n$-dimensional operator space $E$ with $\lambda_{cb}(E) = \sqrt{n}$
or $ex(E) = \sqrt{n}$, at least, at the time of this writing.

In this paper we are going to investigate the extremal cases of $\lambda_{cb}(E) = \sqrt{n}$
and $ex(E) = \sqrt{n}$ and prove the following theorem.

\begin{thm}\label{thm-strict-estimate}
Let $n\geq 2$ and $E \subseteq B(H)$ be an $n$-dimensional
operator space. Then, we have $\lambda_{cb}(E) = \sqrt{n}$ if and
only if $ex(E) = \sqrt{n}$. In other words, $\lambda_{cb}(E) < \sqrt{n}$ if and
only if $ex(E) < \sqrt{n}$.
\end{thm}

$\lambda_{cb}(E)$ is the operator space analogue of the projection
constant $\lambda(X)$ of a Banach space $X$ given by $\lambda(X) =
\sup\{\lambda(X, Y) | X \subseteq Y\}$, where $$\lambda(X, Y) =
\inf\{\norm{P} | P:Y\rightarrow Y \,\,\text{projection onto}\,\, X
\}.$$ See \cite{KL, KT1, KT2} for more information about Banach
space cases and \cite{J2,PS} for operator space cases.

Throughout this paper, we assume that the reader is familiar with
the standard materials about operator spaces (\cite{ER,P1}), completely nuclear maps (\cite{ER}) and completely $p$-summing maps (\cite{P2}).
For a linear map $T : E \rightarrow F$ between operator spaces and
$1\leq p < \infty$ we denote the completely nuclear norm and the completely $p$-summing norm of $T$ by $\nu^o(T)$ and $\pi^o_p(T)$, respectively.

For an index set $I$, $OH(I)$ denote the operator Hilbert space on $\ell_2(I)$
which is introduced in \cite{P1}. When $I = \{1,\cdots, n \}$ for $n\in \mathbb{N}$, we simply write $OH_n$. For a family of operator spaces $(E_i)_{i\in I}$ and an ultrafilter $\U$ on $I$ we denote the ultraproduct of $(E_i)_{i\in I}$ with respect to $\U$ by $\prod_{\U}E_i$.

\section{Proof of the main result}

In the proof we need several lemmas. The first one is about the inclusion between completely
1-summing map and completely 2-summing maps.

\begin{lem}\label{inclusion}
Let $v : E \rightarrow F$ be a completely 1-summing map. Then, $v$
is completely 2-summing with $\pi^o_2(v)\leq \pi^o_1(v).$
\end{lem}
\begin{proof}
Let $E \subseteq B(H)$ for some Hilbert space $H$. Then by Remark
5.7 of \cite{P2} we have an ultrafilter $\U$ over an index set $I$
and the families of positive operators $(a_{\alpha})_{\alpha \in
I}$, $(b_{\alpha})_{\alpha \in I}$ in the unit ball of $S_2(H)$
such that the following diagram commute for some $u$ with
$\norm{u}_{cb} \leq \pi^o_1(v)$ :
\begin{eqnarray}\label{factorization}
\begin{CD}
E @> v >> F \\
@V i VV @AA u A\\
E_{\infty} @>> \mathcal{M} > E_1,
\end{CD}
\end{eqnarray}
where $E_{\infty} = i(E)$ for the complete isometry $$i : B(H)
\hookrightarrow \prod_{\U} B(H), x \mapsto (x)_{\alpha\in I},$$ $E_1
= \overline{Mi(E)}$ (closure in $\prod_{\U} S_1(H)$) for $$M : \prod_{\U}
B(H) \rightarrow \prod_{\U} S_1(H), (x_{\alpha}) \mapsto
(a_{\alpha}x_{\alpha}b_{\alpha})$$ and $\M = M|_{E_{\infty}}$.

Now we split $M = T_2 T_1$, where $$T_1 : \prod_{\U} B(H)
\rightarrow \prod_{\U} S_2(H), (x_{\alpha}) \mapsto
(a^{\frac{1}{2}}_{\alpha}x_{\alpha}b^{\frac{1}{2}}_{\alpha})$$ and
$$T_2 : \prod_{\U} S_2(H) \rightarrow \prod_{\U} S_1(H), (x_{\alpha}) \mapsto
(a^{\frac{1}{2}}_{\alpha}x_{\alpha}b^{\frac{1}{2}}_{\alpha}).$$
Note that
\begin{equation}\label{T2}
\norm{T_2}_{cb} \leq \lim_{\U}\norm{M_{\alpha} : S_2(H)
\rightarrow S_1(H)\,\, , \,\, x \mapsto a^{\frac{1}{2}}_{\alpha}x
b^{\frac{1}{2}}_{\alpha}}_{cb}\leq 1 \end{equation} since
$M_{\alpha}^* = N_{\alpha}$ for $$N_{\alpha} : B(H) \rightarrow
S_2(H)\,\, , \,\, x \mapsto a^{\frac{1}{2}}_{\alpha}x
b^{\frac{1}{2}}_{\alpha}$$ and $\norm{N_{\alpha}}_{cb}\leq 1$.
Thus we have by Theorem 5.1 of \cite{P2} that
\begin{eqnarray*}
\begin{split}
\norm{(vx_{ij})}_{M_n(F)} & = \norm{(uT_2T_1 ix_{ij})}_{M_n(F)}
\leq \pi^o_1(v)\norm{(T_2T_1 ix_{ij})}_{M_n(S_1(H))}\\ & \leq
\pi^o_1(v)\norm{(T_1 ix_{ij})}_{M_n(S_2(H))} =
\pi^o_1(v)\norm{(a_{\alpha}^{\frac{1}{2}}x_{ij}b_{\alpha}^{\frac{1}{2}})}_{M_n(S_2(H))}
\end{split}
\end{eqnarray*}
for any $n\in \mathbb{N}$ and $(x_{ij}) \in M_n(F)$, which implies
$\pi^o_2(v)\leq \pi^o_1(v).$

\end{proof}

The second one is about the trace duality of completely 2-summing norm.
\begin{lem}\label{self-dual}
Let $E$ and $F$ be operator spaces and $E$ be finite dimensional.
Then for $v:F\rightarrow E$ we have
$$(\pi^o_2)^*(v) := \sup\{\abs{{\rm tr}(vu)}|\pi^o_2(u:E\rightarrow F) \leq 1\} = \pi^o_2(v).$$
\end{lem}
\begin{proof} See Lemma 4.7 of \cite{L}
\end{proof}

The last one is about the relationship of the trace and the completely nuclear norm of a linear map on an operator space
with the operator space approximation property.

\begin{lem}\label{trace-com-nuclear}
Let $E$ be an operator space with the operator space approximation property.
Then for any completely nuclear map $u : E\rightarrow E$ we can define ${\rm tr}(u)$, the trace of $u$, and we have
$$\abs{{\rm tr}(u)} \leq \nu^o(u).$$
\end{lem}
\begin{proof} Since $E$ has the operator space approximation property the canonical mapping
$$\Phi : E \widehat{\otimes} E^* \rightarrow E \otimes_{\min}E^*$$ is one-to-one by Theorem 11.2.5 of \cite{ER},
where $\widehat{\otimes}$ (resp. $\otimes_{\min}$) is the projective (resp. injective) tensor product in the category of operator space.
Thus, $\N^o(E)$, the set of all completely nuclear maps on $E$ can be identified with $E \widehat{\otimes}E^*$ with the same norm.
Since we have trace functional define on $E \widehat{\otimes}E^*$ (7.1.12 of \cite{ER}) we can translate it to $\N^o(E)$,
so that we have $$\abs{{\rm tr}(u)} \leq \norm{U}_{E \widehat{\otimes}E^*} = \nu^o(u),$$
where $U \in E \widehat{\otimes}E^*$ is the element associated to $u \in \N^o(E)$.
\end{proof}

Let $E$ and $F$ are operator spaces. Then,
$\Gamma_{\infty}$-norm and $\gamma_{\infty}$-norm of a linear map
$v: E\rightarrow F$ are defined by
$$\Gamma_{\infty}(v) = \inf\norm{\alpha}_{cb}\norm{\beta}_{cb},$$
where the infimum is taken over all Hilbert space $H$ and the factorization $$v : E \stackrel{\alpha}{\rightarrow}
B(H) \stackrel{\beta}{\rightarrow} F,$$ and
$$\gamma_{\infty}(v) = \inf\norm{\alpha}_{cb}\norm{\beta}_{cb},$$
where the infimum is taken over all $m\in \mathbb{N}$ and the factorization $$v : E \stackrel{\alpha}{\rightarrow}
M_m \stackrel{\beta}{\rightarrow} F.$$
See section 4 of \cite{J2} or \cite{EJR} for the details.

Now we are ready to prove our main result. The proof follows the classical idea of
\cite{KL}.

\vspace{0.5cm}
{\it proof of Theorem \ref{thm-strict-estimate} :}
\vspace{0.5cm}

By Fact \ref{fact1} and Fact \ref{fact2} it is enough to show that the condition $\lambda_{cb}(E) =
\sqrt{n}$ is contradictory to the condition $ex(E) = d_{\mathcal{S}\K}(E) < \sqrt{n}$.

\vspace{0.5cm}

{\it Step 1 :} $\pi^o_1(I_E) = \sqrt{n}.$

\vspace{0.5cm}

By trace duality and Lemma 4.1 and 4.2 of \cite{J2} (or see
Theorem 7.6 of \cite{EJR}) we have
$$\lambda_{cb}(E) = \Gamma_{\infty}(I_E) = \gamma_{\infty}(I_E) = \sup_{u\in \pi^o_1(E)}
\frac{\abs{\text{tr}(u)}}{\pi^o_1(u)}.$$ Since $E$ is finite
dimensional, we can find $u\in CB(E)$ such that $$
\frac{\abs{\text{tr}(u)}}{\pi^o_1(u)} = \sqrt{n}$$ and by multiplying
suitable constant we can also assume that $\pi^o_2(u) = \sqrt{n}$.
Then we have by Lemma \ref{inclusion}, Lemma \ref{self-dual} and
Theorem 6.13 of \cite{P2} that
$$ n  = \sqrt{n} \pi^o_2(u) \leq \sqrt{n}\pi^o_1(u) = \abs{\text{tr}(u)}
\leq \pi^o_2(u)\pi^o_2(I_E) = n.$$
Thus, we get $$\pi^o_1(u) = \sqrt{n}\,\, \text{and}\,\, \abs{\text{tr}(u)} = n.$$ Now
we will show that $u$ is actually $I_E$. By Proposition 6.1 of
\cite{P2} we have the factorization $$v : E \stackrel{A}{\rightarrow} OH_n \stackrel{B}{\rightarrow}
E\,\, \text{with} \,\,\pi^o_2(A)\norm{B}_{cb} \leq
\sqrt{n}.$$ If we let $v : OH_n \rightarrow OH_n$ by $v = AB$, we
have $\text{tr}(v) = \text{tr}(v^*) = \text{tr}(u)$ and
\begin{eqnarray*}
\begin{split}
\norm{I_{OH_n}-v}^2_{HS} & =
\text{tr}\big((I_{OH_n}-v)(I_{OH_n}-v)^*\big)\\
& = \text{tr}(I_{OH_n}) -2\text{tr}(u) + \text{tr}(vv^*)\\
& = n -2n + \norm{v}^2_{HS} = (\pi^o_2(v))^2 -n\\ & \leq
(\pi^o_2(A)\norm{B}_{cb})^2 -n \leq 0,
\end{split}
\end{eqnarray*}
which leads us to our desired conclusion.

\vspace{0.5cm}

{\it Step 2 :} Now we factorize $I_E$ as in the proof of Lemma
\ref{inclusion}. Then we have an ultrafilter $\U$, the families of
positive operators $(a_{\alpha})_{\alpha \in I}$,
$(b_{\alpha})_{\alpha \in I}$ in the unit ball of $S_2(H)$ such
that the diagram (\ref{factorization}) commute for some $u$ with
$$\norm{u}_{cb} \leq \pi^o_1(I_E) = \sqrt{n}.$$ Then we can find a
rank $n$ projection $$w_1 : i(B(H)) \rightarrow i(B(H))\,\, \text{onto}\,\, E_{\infty}\,\, \text{with}\,\,
\pi^o_1(w_1) \leq \sqrt{n}.$$

\vspace{0.5cm}

Consider $iu : E_1 \rightarrow i(B(H))$. Since $i$ is a complete
isometry, $i(B(H))$ is injective in the operator space sense, so
that we can extend $iu$ to $$\tilde{u} : \prod_{\U} S_1(H)
\rightarrow i(B(H))\,\, \text{with}\,\, \norm{\tilde{u}}_{cb} =
\norm{iu}_{cb}.$$ Now we consider the same factorization $M = T_2
T_1$ as before. Note that $$\pi^o_2(T_1) \leq 1 \,\, \text{and}\,\,
\norm{T_2}_{cb} \leq 1$$ by the same
calculation as the proof for (5.8) of \cite{P2} and (\ref{T2}), respectively. Then for $$w := T_1
\tilde{u}T_2 : \prod_{\U} S_2(H) \rightarrow \prod_{\U} S_2(H)$$ we
have
\begin{equation}\label{w-HS-norm}
\norm{w}_{HS} = \pi^o_2(w)  \leq
\pi^o_2(T_1)\norm{\tilde{u}}_{cb}\norm{T_2}_{cb} \leq \pi^o_2(T_1)\norm{u}_{cb} \leq \sqrt{n}.
\end{equation}
Since $T_1 i$ is 1-1, $F := T_1 i(E)$ is $n$-dimensional.
Furthermore, since $$wT_1ix = T_1\tilde{u}T_2 T_1 ix = T_1 iu\M ix
= T_1 ix$$ for all $x\in E$, we have $w|_{F} = I_F$, which means
$\abs{\lambda_k(w)} \geq 1$ for $1 \leq k \leq n$,
where $(\lambda_k(w))_{k\geq 1}$ is the sequence of eigenvalues of $w$ which is non-increasing in absolute value
and repeated as often as its multiplicity. By applying
Weyl's inequality (Lemma 3.5.4 of \cite{Pie}) and (\ref{w-HS-norm}), we get
$$n \leq \sum^n_{k=1}\abs{\lambda_k(w)}^2 \leq
\sum^{\infty}_{k=1}s_k(w)^2 = \norm{w}^2_{HS} \leq n,$$ where $(s_k(w))_{k\geq 1}$ is the sequence of singular values of $w$.
Then we have
$$\abs{\lambda_k(w)} =
\begin{cases} \,1 \text{\,\,\,\,\,\,\,\,\,\,if $1 \leq k \leq n$,}\\
\,0 \text{\,\,\,\,\,\,\,\, if $k > n$,}
\end{cases}$$
which implies $w$ has rank at most $n$ and so does $$w_1 :=
\tilde{u}\M = \tilde{u}T_2 T_1|_{i(B(H))} : i(B(H)) \rightarrow
i(B(H)).$$ Actually, $w_1$ is our desired rank $n$ projection.
Indeed, we have
$$w_1ix = \tilde{u}\M ix = iu\M ix = ix$$ for all $x \in E$ and
since $E_{\infty}$ is $n$-dimensional, $w_1$ maps onto
$E_{\infty}$. Moreover, we have $$\pi^o(w_1) \leq
\norm{\tilde{u}}_{cb}\pi^o_1(\M) \leq \sqrt{n}$$ since
$\pi^o_1(\M) \leq 1$. ((5.7) of \cite{P2})

\vspace{0.5cm}

{\it Step 3 :} Since $d_{\mathcal{S}\K}(E_{\infty}) =
d_{\mathcal{S}\K}(E) < \sqrt{n}$, we have $F \in \K$ and an
isomorphism $$T : E_{\infty} \rightarrow F\,\, \text{with}\,\,
\norm{T}_{cb}\norm{T^{-1}}_{cb} < \sqrt{n}.$$ By the fundamental
extension theorem (Theorem 1.6 of \cite{P3}) we have extensions
$$\widetilde{T}: i(B(H)) \rightarrow B(\ell_2)\,\,\, \text{and}\,\,\,
\widetilde{T^{-1}}: B(\ell_2) \rightarrow i(B(H))$$ of $T$ and
$T^{-1}$, respectively, with $\norm{\widetilde{T}}_{cb} =
\norm{T}_{cb}$ and $\norm{\widetilde{T^{-1}}}_{cb} =
\norm{T^{-1}}_{cb}$.

Let $\tilde{w_1} = \widetilde{T}w_1\widetilde{T^{-1}} : B(\ell_2)
\rightarrow B(\ell_2)$. Then clearly we have
$\text{ran}(\tilde{w_1}) \subseteq F$ and $\tilde{w_1}|_F = I_F$,
which means that $\tilde{w_1}$ is also a rank $n$ projection from
$B(\ell_2)$ onto $F$. Since $F \subseteq \K$ and $\K$ satisfies the operator space approximation property
we have by Lemma \ref{trace-com-nuclear} and Corollary 15.5.4 of \cite{ER} that
\begin{eqnarray*}
\begin{split}
n = \abs{\text{tr}(\tilde{w_1}|_{\K} : \K \rightarrow \K)} & \leq
\nu^o(\tilde{w_1}|_{\K}: \K \rightarrow \K) =
\pi^o_1(\tilde{w_1}|_{\K}: \K \rightarrow \K)\\ & =
\pi^o_1(\tilde{w_1}|_{\K}: \K \rightarrow B(\ell_2)) \leq
\norm{\widetilde{T}}_{cb}\norm{\widetilde{T^{-1}}}_{cb}\pi^o_1(w_1)\\
& \leq \norm{T}_{cb}\norm{T^{-1}}_{cb}\sqrt{n} < n,
\end{split}
\end{eqnarray*}
which is a contradiction.

\bibliographystyle{amsplain}
\providecommand{\bysame}{\leavevmode\hbox
to3em{\hrulefill}\thinspace}

\end{document}